\documentclass[smallextended,numbook,runningheads]{svjour3}     
\smartqed  
\usepackage{graphicx}
\usepackage{amsmath}
\usepackage{epstopdf} 
\usepackage{mathptmx}      
\usepackage{amsfonts}
\DeclareMathOperator{\cay}{cay}
\DeclareMathOperator{\ad}{ad}
\DeclareMathOperator{\diag}{diag}

\newtheorem{lgrthm}[theorem]{Algorithm}
\journalname{BIT}
\begin{document}

\title{Higher Strong Order Methods for Itô SDEs on Matrix Lie Groups}


\author{Michelle Muniz \and Matthias Ehrhardt \and Michael Günther \and Renate Winkler}


\institute{Chair of Applied Mathematics / Numerical Analysis,\\ 
Faculty of Mathematics and Natural Sciences, \\
Bergische Universität Wuppertal,\\ 
Gaußstrasse 20, 42119 Wuppertal, Germany \\
\email{\{muniz,ehrhardt,guenther,winkler\}@uni-wuppertal.de}
}

\date{Received: date / Accepted: date}

\maketitle

\begin{abstract}
In this paper we present a general procedure for designing higher strong order methods for Itô stochastic differential equations on matrix Lie groups and illustrate this strategy with two novel schemes that have a strong convergence order of 1.5. 
 Based on the Runge-Kutta--Munthe-Kaas (RKMK) method for ordinary differential equations on Lie groups, we present a stochastic version of this scheme and derive a condition such that the stochastic RKMK has the same strong convergence order as the underlying stochastic Runge-Kutta method.
 Further, we show how our higher order schemes can be applied in a mechanical engineering as well as in a financial mathematics setting.
\end{abstract}

\section{Introduction}

In recent years, more and more interrelations in mechanics and finance have been modeled by stochastic differential equations (SDEs) on Lie groups. 
A trend can be observed that shows that kinematic models, which were previously expressed by ordinary differential equations (ODEs), are now extended by terms that include stochastic processes in order to include possible stochastic perturbations. 
Examples can be found in the modeling of rigid bodies like satellites, vehicles and robots \cite{Chir1,Chir12,MaWi08,ZuVa19}. 
Furthermore, SDEs on Lie groups are also considered in the estimation of object motion from a sequence of projections~\cite{Soat94} 
and in the representation of the precessional motion of magnetization in a solid~\cite{AbBu17}. 

In financial mathematics, the consideration of stochastic processes is essential, and the solution of SDEs has been performed for many years, but usually not on Lie groups. 
However, the use of Lie groups to solve existing or to create new financial models could be of central importance for dealing with geometric constraints. We are confronted with geometric constraints, e.g.\ in the form of a positivity constraint on interest rates~\cite{Park10,LiPr16,Wang20} or a symmetry and positivity constraint on covariance and correlation matrices~\cite{Muniz20}, which are important e.g.\ in risk management and portfolio optimization.

Despite these diversified applications, the available literature on analysis and numerical methods for SDEs on Lie groups is limited, in contrast to the available literature on ODEs on Lie groups (e.g.\ \cite{Munt98,Munt99,CrGr93,Hair06,Iser05,Cell14}). 
Furthermore, the available literature on Lie group SDEs mainly concerns Stratonovich SDEs \cite{BuBu99,MaWi08,AbBu17,Wang20}. 
Readers interested in Itô SDEs on Lie groups will only find the \textit{geometric Euler-Maruyama} scheme with strong order $\gamma=1$ appearing in \cite{Marj15,MaSo18,PiSo16} and more recent the existence and convergence proof of the stochastic Magnus expansion in \cite{Kamm20}.
However, the consideration of Itô SDEs is crucial for its application in finance and due to the geometric constraints Stratonovich SDEs on matrix Lie groups cannot simply be transformed into Itô SDEs as in the traditional, non-geometric case.

Our contribution to this field of research is a general procedure on how to set up structure-preserving schemes of higher strong order for Itô SDEs on matrix Lie groups. Based on the \textit{Magnus expansion} we apply Itô-Taylor schemes or stochastic Runge-Kutta (SRK) schemes to solve a corresponding SDE in the Lie algebra. 
Using a SRK method can be interpreted as a stochastic version of Runge-Kutta--Munthe-Kaas (RKMK) methods. 
Under these circumstances, we derive a condition such that the stochastic RKMK scheme inherits the strong convergence order $\gamma$ of the SRK method applied in the Lie algebra.

The remainder of the paper is organized as follows.
We start with an introduction to matrix Lie groups, their corresponding Lie algebras and the linear Itô matrix SDE which we consider in this geometric setting in Section~\ref{sec2}.
In Section~\ref{sec3} we take a closer look on how SDEs on Lie groups can be solved numerically and present our higher strong order methods. 
Then we provide some numerical and application examples in Section~\ref{sec4}. A conclusion of our results is given in Section~\ref{sec5}.

\section{SDEs on Matrix Lie Groups}\label{sec2}
A Lie group is a differentiable manifold, which is also a group $G$ with a differentiable product that maps $G\times G\to G$. 
Matrix Lie groups are Lie groups, which are also subgroups of GL($n$) for $n\in\mathbb{N}$. 
The tangent space at the identity of a matrix Lie group $G$ is called Lie algebra $\mathfrak{g}$. 
The Lie algebra is closed under forming Lie brackets $[\cdot,\cdot]$ (also called commutators) of its elements.
For further details on Lie groups and Lie algebras we refer the interested reader to \cite{Hall15}.

On a matrix Lie group $G$ we consider the linear matrix-valued Itô SDE
\begin{equation}\label{Q SDE}
    dQ_t = Q_t K_t\,dt + Q_t  V_{t}\,dW_{t}, \quad Q_0=I_{n\times n},
\end{equation}
where $K_t, V_{t} \in \mathbb{R}^{n \times n}$ are given coefficient matrices, 
$W_t$ denotes the standard Brownian motion, i.e.\ it holds $dW_t\sim\mathcal{N}(0,dt)$ 
and $I_{n\times n}$ is the $n$-dimensional identity matrix.
In general, there exists no closed form solution to \eqref{Q SDE}.
However, a solution can be defined via a \textit{Magnus expansion} 
$Q_t=Q_0\,\psi(\Omega_t)$ (see \cite{Hair06,Magn54,MaSo18}), 
where $\Omega_t\in\mathbb{R}^{n\times n}$ obeys the following matrix SDE
\begin{equation}\label{Omega SDE}
	d\Omega_t = A(\Omega_t)\, dt + \Gamma(\Omega_t) \,dW_t, \quad \Omega_0 = 0_{n\times n}.
\end{equation}
The drift and diffusion coefficient are given by 
\begin{equation}\label{eq: Omega coefficients}
   A(\Omega_t) = d\psi_{-\Omega_t}^{-1}
   \Bigl(K_t - \frac{1}{2}V_t^2 - \frac{1}{2}C(\Omega_t)\Bigr),\quad \Gamma(\Omega_t) 
   = d\psi_{-\Omega_t}^{-1}(V_t)
\end{equation}
with 
\begin{equation}\label{eq: C coefficient}
   C(\Omega_t) = \Bigl(\frac{d}{d\Omega_t}d\psi_{-\Omega_t}\bigl(\Gamma(\Omega_t)\bigr)\Bigr)\,\Gamma(\Omega_t)
\end{equation}
(see appendix for proof). 
For $Q_t \in G$, the solution of the matrix SDE~\eqref{Omega SDE} $\Omega_t$ 
is an element of the Lie algebra $\mathfrak{g}$. 
The mapping $\psi\colon\mathfrak{g}\to G$ is considered to be a local 
diffeomorphism between the Lie algebra and the corresponding Lie group near $\Omega=0_{n\times n}$. 

\subsection{The exponential map as local parametrization}
A common choice for $\psi(\Omega)$ is $\exp(\Omega)=\sum_{k=0}^{\infty}\frac{1}{k!}\Omega^k$ with the derivative
\begin{equation*}
    \Bigl(\frac{d}{d\Omega} \exp(\Omega)\Bigr)H 
    = \bigl(d\exp_{\Omega}(H)\bigr)\exp(\Omega)
    =\exp(\Omega)\left(d\exp_{-\Omega}(H)\right)
\end{equation*}
where
\begin{equation}\label{eq: dexp}
    d\exp_{-\Omega}(H) = \sum_{k=0}^{\infty}\frac{1}{(k+1)!}\ad_{-\Omega}^k(H).
\end{equation}
The inverse of $d\exp$ is given in the following Lemma \cite[p.~84]{Hair06}.
\begin{lemma}[Baker, 1905]
    If the eigenvalues of the linear operator $\ad_{\Omega}$ are different from $2\ell\pi i$ with $\ell \in \{\pm1,\pm2,\dots\}$, 
    then $d\exp_{-\Omega}$ is invertible. Furthermore, we have for $\|\Omega\|<\pi$ that
    \begin{equation}\label{eq: dexpinv}
        d\exp_{-\Omega}^{-1}(H) = \sum_{k=0}^{\infty} \frac{B_k}{k!}\ad_{-\Omega}^k(H),
    \end{equation}
    where $B_k$ are the Bernoulli numbers, defined by $\sum_{k=0}^{\infty}(B_k/k!)x^k=x/(e^x-1)$.
\end{lemma}
We recall that the first three Bernoulli numbers
are given by $B_0=1$, $B_1=-\frac{1}{2}$, $B_2=\frac{1}{6}$ and that $B_{2m+1}=0$ holds for $m\in \mathbb{N}$. 

By $\ad_{\Omega}(H)=\lbrack \Omega,H\rbrack = \Omega H - H\Omega$ we express
the adjoint operator which is used iteratively
	\begin{equation*}
    	\ad_{\Omega}^0(H) =H,\quad
    	\ad_{\Omega}^k(H)=\bigl\lbrack \Omega, \ad_{\Omega}^{k-1}(H)\bigr\rbrack
    	=\ad_{\Omega}\bigl(\ad_{\Omega}^{k-1}(H)\bigr), \quad k\ge 1. 
	\end{equation*}

With these expressions and Itô rules the coefficient in \eqref{eq: C coefficient} can be simplified to
\begin{equation}\label{eq: C exp}
    C(\Omega_t) = \sum_{p=0}^{\infty} \sum_{q=0}^{\infty} \frac{1}{(p+q+2)}\frac{(-1)^p}{p!(q+1)!}
    \ad_{\Omega_t}^p\Bigl(\ad_{\Gamma(\Omega_t)}\bigl(\ad_{\Omega_t}^q\bigl(\Gamma(\Omega_t)\bigr)\bigr)\Bigr),
\end{equation}
we refer to \cite{MaSo18} for the concise derivation of this expression.

\subsection{The Cayley map as local parametrization}
With all these series related to $\psi=\exp$, the question arises whether there is another mapping $\psi\colon\mathfrak{g} \to G$, 
which is not based on the evaluation of an infinite number of summands. 
In case of a \textit{quadratic Lie group} the answer is yes, there is a mapping, 
namely the \textit{Cayley transformation} 
\begin{equation*}
    \cay(\Omega) = (I - \Omega)^{-1}(I + \Omega).
\end{equation*}
A quadratic Lie group $G$ is a set of matrices $Q$ 
that fulfill the equation $Q^{\top}PQ = P$ for a given constant matrix $P$. 
For the derivative of $\cay(\Omega)$ we have
\begin{equation*}
    \Bigl(\frac{d}{d\Omega} \cay(\Omega)\Bigr)H = \bigl(d\cay_{\Omega}(H)\bigr)\cay(\Omega)
    =\cay(\Omega)\left(d\cay_{-\Omega}(H)\right).
\end{equation*}
The analogue expression to \eqref{eq: dexp} reads
\begin{equation*}
    d\cay_{-\Omega}(H) = 2(I + \Omega)^{-1} H (I-\Omega)^{-1}
\end{equation*}
with the inverse given by
\begin{equation}\label{eq: dcayinv}
    d\cay_{-\Omega}^{-1}(H) = \frac{1}{2}(I + \Omega) H (I-\Omega),
\end{equation}
see \cite{Hair06}. Using the Cayley map as local parametrization the coefficient 
$C(\Omega_t)$ in \eqref{eq: C coefficient} is given by
\begin{equation}\label{eq: C cay}
    C(\Omega_t) = V_t\Omega_t V_t
\end{equation}
(see appendix for proof).

\subsection{Example: SDEs on SO($n$)}
As an example for a matrix Lie group, we take a closer look on the special orthogonal group
\begin{equation*}
	\mathrm{SO}(n)=\{X \in \mathrm{GL}(n) : X^{\top}X=I, \, \det(X)=1\}
\end{equation*}
which is a quadratic Lie group, such that the Cayley map is also applicable as a local parametrization. 
The corresponding Lie algebra consists of skew-symmetric matrices,
\begin{equation*}
    \mathfrak{so}(n) = \{Y \in \mathrm{GL}(n) : Y + Y^\top = 0\}.
\end{equation*}
Since we are interested in structure preservation we need conditions that tell us when the solution of an SDE on SO($n$) is kept on the manifold.

\begin{theorem}\label{theorem: so(n)}
For the solution $Q_t$ of \eqref{Q SDE} it holds $Q_t \in$ SO($n$) if and only if the coefficient matrices satisfy $V_t \in \mathfrak{so}(n)$ and $K_t+K_t^{\top}=V_t^2$.
\end{theorem}
For the proof of this theorem we refer to~\cite{MaSo18}.

\section{Numerical methods for SDEs on Lie Groups}\label{sec3}

Applying standard numerical methods for SDEs directly to the linear matrix-valued 
Itô SDE~\eqref{Q SDE} will result in a \textit{drift off}, 
i.e.\ the numerical approximations do not stay on the manifold. 
Consequently, one needs to consider special numerical methods 
that preserve the geometric properties of the Lie group $G$. 

As the Lie algebra $\mathfrak{g}$ represents a linear space with Euclidean-like geometry, it appears reasonable to compute the numerical approximations of the matrix SDE~\eqref{Omega SDE} and to project the solution back onto the Lie group $G$.

A simple scheme based on the Runge-Kutta--Munthe-Kaas schemes for ODEs \cite{Munt98} that puts the 
described approach into practice can be found in~\cite{Marj15} and is presented in the following algorithm.

\begin{lgrthm}\label{algo}
	Divide the time interval $\lbrack 0,T\rbrack$ uniformly into $J$ subintervals 
	$\lbrack t_j, t_{j+1} \rbrack$, $j=0,1,\dots,J-1$ and define the time step $\Delta = t_{j+1} - t_j$. 
	Let $Q_t=Q_0\psi(\Omega_t)$ with $\psi\colon\mathfrak{g}\to G$ be a local parametrization of the Lie group $G$.
	Starting with $t_0=0$, $Q_0=I_{n\times n}$ and $\Omega_0=0_{n\times n}$ the 
	following steps are repeated over successive intervals $\lbrack t_j, t_{j+1}\rbrack$ until $t_{j+1}=T$.
	\begin{enumerate}
		\item \textbf{Initialization step:} Let $Q_j$ be the approximation of $Q_t$ at time $t=t_j$.
		\item \textbf{Numerical method step:} Compute an approximation $\Omega_1\approx\Omega_{\Delta}$ by applying a stochastic 
		Itô-Taylor or stochastic Runge-Kutta method to the matrix SDE~\eqref{Omega SDE}. 
		\item \textbf{Projection step:} Set $Q_{j+1}=Q_j\, \psi(\Omega_1)$.
	\end{enumerate}
\end{lgrthm}

The order of convergence of these Lie group structure-preserving schemes clearly depends 
on the numerical method used in the second step of the algorithm. 
In order to analyze the accuracy of our geometric numerical methods we recall that an approximating 
process $X_t^{\Delta}$ is said to \textit{converge in a strong sense with order} $\gamma > 0$ to the 
Itô process $X_t$ if there exists a finite constant $K$ and a $\Delta'>0$ such that
\begin{equation}\label{eq: order}
    \mathbb{E}[|X_T - X_T^{\Delta}|] \leq K \Delta^{\gamma}
\end{equation}
for any time discretization with maximum step size $\Delta \in (0,\Delta')$ \cite{KlPl92}.

\subsection{Geometric schemes of strong order 1}\label{subsec: order1}
Using the Euler-Maruyama scheme in the numerical method step of Algorithm~\ref{algo} results in
\begin{align}
    \begin{split}\label{eq: gEM}
        \Omega_{1} &= \Omega_0 + A(\Omega_0) \Delta + \Gamma(\Omega_0) \Delta W \\
        &= d\psi_{-\Omega_0}^{-1}\Bigl(K_j - \frac{1}{2}V_j^2 \Bigr) \Delta +  d\psi_{-\Omega_0}^{-1}(V_j) \Delta W, \\
        Q_{j+1} &= Q_j \psi(\Omega_1),
    \end{split}
\end{align}
where $\Delta W \sim \mathcal{N}(0,\Delta)$. Note that $C(\Omega_0)=0_{n\times n}$ for both mappings $\psi=\exp$ (see \eqref{eq: C exp}) and $\psi=\cay$ (see \eqref{eq: C cay}) which is why we neglect this coefficient from here on.

Since this scheme \eqref{eq: gEM} preserves the geometry of the Lie group $G$ it was called the \textit{geometric Euler-Maruyama scheme}~\cite{MaSo18}.
It can be specified according to the mapping.

For $\psi=\exp$, we get
\begin{align}
    \begin{split}\label{eq: gEM exp}
        \Omega_{1} &= d\exp_{-\Omega_0}^{-1}\Bigl(K_j - \frac{1}{2}V_j^2 \Bigr) \Delta +  d\exp_{-\Omega_0}^{-1}(V_j) \Delta W \\
        &= \Bigl(K_j - \frac{1}{2}V_j^2 \Bigr) \Delta +  V_j \Delta W, \\
        Q_{j+1} &= Q_j\,\exp(\Omega_1),
    \end{split}
\end{align}
where inserting $\Omega_0=0_{n \times n}$ is equivalent to truncating the infinite 
series~\eqref{eq: dexpinv} after the first summand,
right before any dependence on $\Omega$ appears.

Using $\psi=\cay$ instead, we obtain
\begin{align*}
\begin{split}
    \Omega_1 &= d\cay_{-\Omega_0}^{-1}\Bigl(K_j-\frac{1}{2} V_{j}^2\Bigr) \Delta +  d\cay_{-\Omega_0}^{-1}(V_{j})\,\Delta W \\
    &= \frac{1}{2}\Bigl(K_j-\frac{1}{2} V_{j}^2\Bigr)\Delta + \frac{1}{2}V_j \Delta W, \\
    Q_{j+1} &= Q_j\,\cay(\Omega_1) = Q_j(I - \Omega_1)^{-1}(I+ \Omega_1).
\end{split}
\end{align*}

In both cases we see that the diffusion term is only dependent on time and not on the solution itself. This is called \textit{additive noise}~\cite{KlPl92} and it is the reason why these schemes have strong order $\gamma=1$ instead of $\gamma=0.5$ as expected for the traditional Euler-Maruyama method. A general proof of the geometric Euler-Maruyama method converging with strong order $\gamma=1$ can be found in~\cite{PiSo16}.

\subsection{Geometric schemes of higher order}
A higher strong order than $\gamma = 1$ can be achieved by applying e.g.~the strong 
Itô-Taylor approximation of order $\gamma=1.5$ (see \cite{KlPl92}) in the second step of Algorithm~\ref{algo}. 
By doing so, we obtain
\begin{equation}
\begin{split}\label{eq: ItôTaylor1.5}
    \Omega_1 =&~ A(\Omega_0) \Delta + \Gamma(\Omega_0) \,\Delta W + \frac{1}{2}\Gamma'\Gamma(\Omega_0)\bigl( (\Delta W)^2 - \Delta \bigr)\\
    &+ A'\Gamma(\Omega_0)\Delta Z \\
    &+ \frac{1}{2}\left(A'A(\Omega_0) + \frac{1}{2}A''\Gamma^2(\Omega_0)\right)\Delta^2 \\
    &+ \left(\Gamma'A(\Omega_0) + \frac{1}{2}\Gamma''\Gamma^2(\Omega_0)\right)(\Delta W \Delta - \Delta Z) \\
    &+ \frac{1}{2}\left(\Gamma'\Gamma(\Omega_0)\right)'\Gamma(\Omega_0)\; \Bigl(\frac{1}{3}(\Delta W)^2 - \Delta\Bigr)\,\Delta W,\\
    Q_{j+1} =&~ Q_j \,\psi(\Omega_1).
\end{split}
\end{equation}
Representing the double integral $\int_{\tau_{\ell}}^{\tau_{\ell+1}} \int_{\tau_{\ell}}^{s_2} dW_{s_1}ds_2$, 
the random variable $\Delta Z$ is normally distributed with mean $\mathbb{E}[\Delta Z]=0$, 
variance $\mathbb{E}\bigl[(\Delta Z)^2\bigr]=\frac{1}{3}\Delta^3$ and covariance $\mathbb{E}[\Delta Z \Delta W]=\frac{1}{2}\Delta^2$.
We consider the matrix derivatives as directional derivatives, e.g.
\begin{equation*}
    A'H = \left(\frac{d}{d\Omega}A(\Omega)\right)H = \frac{d}{d\epsilon}\left.A(\Omega+\epsilon H)\right|_{\epsilon=0}
\end{equation*}
which we then evaluate at $\Omega_0$.
The computation of the needed matrix derivatives for $\psi=\exp$ and $\psi=\cay$ is provided in the Appendix.

A strong order of $\gamma=1.5$ can also be achieved by applying a stochastic Runge-Kutta method of that order to the SDE~\eqref{Omega SDE}. By using the stochastic Runge-Kutta scheme of order $\gamma=1.5$ of R\"oßler \cite{Roes05}, we can avoid computing the derivatives in~\eqref{eq: ItôTaylor1.5} and we obtain
\begin{align}
\begin{split}\label{eq: Roessler1.5}
    \Omega_{1} =&~
    + \left( \frac{1}{3}A(H_1) + \frac{2}{3}A(H_2)\right)\Delta \\
    &+ \left(\frac{13}{4}\Gamma(\Tilde{H}_1) - \frac{9}{4}\Gamma(\Tilde{H}_2) - \frac{9}{4}\Gamma(\Tilde{H}_3) + \frac{9}{4}\Gamma(\Tilde{H}_4)\right)\Delta W \\
    &+ \left(-\frac{15}{4}\Gamma(\Tilde{H}_1) + \frac{15}{4}\Gamma(\Tilde{H}_2) + \frac{3}{4}\Gamma(\Tilde{H}_3) - \frac{3}{4}\Gamma(\Tilde{H}_4)\right)\frac{1}{2\sqrt{\Delta}}\bigl((\Delta W)^2-\Delta\bigr)\\
    &+ \left(-\frac{9}{4}\Gamma(\Tilde{H}_1) + \frac{9}{4}\Gamma(\Tilde{H}_2) + \frac{9}{4}\Gamma(\Tilde{H}_3) - \frac{9}{4}\Gamma(\Tilde{H}_4)\right)\frac{\Delta Z}{\Delta} \\
    &+ \left(6\Gamma(\Tilde{H}_1) - 9\Gamma(\Tilde{H}_2) + 3\Gamma(\Tilde{H}_4)\right)
    \frac{1}{3!\Delta}\bigl((\Delta W)^2-3\Delta\bigr)\Delta W, \\
    Q_{j+1} =&~ Q_j \,\psi(\Omega_1),
\end{split}
\end{align}
with the stage values
\begin{align*}
    H_1 &= H_3 = \Tilde{H}_1 = \Omega_0, \quad H_2 = \frac{3}{4}A(H_1)\Delta + \frac{3}{2}\Gamma(\Tilde{H}_1)\frac{\Delta Z}{\Delta}, \\
    \Tilde{H}_2 &= \frac{1}{9}A(H_1)\Delta + \frac{1}{3}\Gamma(\Tilde{H}_1)\sqrt{\Delta}, \\
    \Tilde{H}_3 &= \frac{5}{9}A(H_1)\Delta + \frac{1}{3}A(H_2)\Delta - \frac{1}{3}\Gamma(\Tilde{H}_1)\sqrt{\Delta} + \Gamma(\Tilde{H}_2)\sqrt{\Delta}, \\
    \Tilde{H}_4 &= A(H_1)\Delta + \frac{1}{3}A(H_2)\Delta + A(H_3)\Delta \\
    &\qquad + \Gamma(\Tilde{H}_1)\sqrt{\Delta} - \Gamma(\Tilde{H}_2)\sqrt{\Delta} + \Gamma(\Tilde{H}_3)\sqrt{\Delta}.
\end{align*}
The exploitation of stochastic Runge-Kutta methods gives us the benefit of a derivative-free scheme. 
However, using the mapping $\psi=\exp$ raises the question of how large the truncation index $q$ must be chosen in the truncated approximation for \eqref{eq: dexpinv},
\begin{equation}\label{eq: dexpinv trunc}
	\sum_{k=0}^{q}\frac{B_k}{k!}\ad_{-\Omega}^k(H) 
	= H - \frac{1}{2} \lbrack -\Omega, H \rbrack + \frac{1}{12} \bigl\lbrack - \Omega, \lbrack -\Omega,H \rbrack \bigr\rbrack + \ldots,
\end{equation}
in order to maintain a strong order of $\gamma=1.5$. 
More generally, a condition is needed which connects the truncation index $q$ with the aimed strong convergence order $\gamma$.

Inspired by \cite[Theorem~IV.8.5.]{Hair06} for Runge-Kutta--Munthe-Kaas methods to solve deterministic matrix ODEs we formulate the following theorem.

\begin{theorem}\label{theorem: index}
Consider Algorithm~\ref{algo} with $\psi=\exp$. 
Let the applied stochastic Runge-Kutta method in the second step of Algorithm~\ref{algo} be of strong order $\gamma$.
If the truncation index $q$ in~\eqref{eq: dexpinv trunc} satisfies $q\ge2\gamma -2$, then the method of Algorithm~\ref{algo} is of strong order $\gamma$.
\end{theorem}

\begin{proof}
  According to the definition of strong convergence~\eqref{eq: order} we have to show that
  \begin{equation*}
      \mathbb{E}[\|\Omega_{\Delta}-\Omega_1\|] \le K\Delta^{(q+2)/2}
  \end{equation*}
  where $\Omega_{\Delta}$ is the exact solution of~\eqref{Omega SDE} with $\psi=\exp$ at $t=\Delta$, $\Omega_1$ is the numerical approximation obtained in the second step of Algorithm~\ref{algo} and $K$ is a finite constant.
  
  Let $\Omega_{\Delta}^q$ be the exact solution of the truncated version of~\eqref{Omega SDE} with $\psi=\exp$ at $t=\Delta$, namely
  \begin{equation*}
      d\Omega_t = \sum_{k=0}^q \frac{B_k}{k!}\ad_{-\Omega_t}^k(K_t-\frac{1}{2}V_t^2) + \sum_{k=0}^q \frac{B_k}{k!}\ad_{-\Omega_t}^k(V_t).
  \end{equation*}
  
  Our proof is divided into six steps.
  
  \subsubsection*{Step 1: Numerical error}
  We consider the absolute error in the Frobenius norm and estimate the error in the $L^1$-norm by the $L^2$-norm. Then, we use the Minkowski inequality by introducing $\Omega_{\Delta}^q$.
  \begin{align*}
      \mathbb{E}[\|\Omega_{\Delta}-\Omega_1\|] 
      &\leq \left(\mathbb{E}\left[\|\Omega_{\Delta}-\Omega_1\|^2\right]\right)^{1/2} \\
      &\leq \left(\mathbb{E}\left[\|\Omega_{\Delta}-\Omega_{\Delta}^q\|^2\right]\right)^{1/2} + \left(\mathbb{E}\left[\|\Omega_{\Delta}^q-\Omega_1\|^2\right]\right)^{1/2}
  \end{align*}
  We are left with the modelling error, which corresponds to the first summand, and the numerical error, the second summand. The numerical error can be estimated by 
  \begin{equation*}
      \left(\mathbb{E}\left[\|\Omega_{\Delta}^q-\Omega_1\|^2\right]\right)^{1/2} \leq \Tilde{K}\Delta^{\gamma} \quad \text{for} \; \Tilde{K}<\infty,
  \end{equation*}
  because we assume that we are applying a SRK method of strong order $\gamma$.
  
  In other words, it remains to be shown that
  \begin{equation*}
      \left(\mathbb{E}\left[\|\Omega_{\Delta}-\Omega_{\Delta}^q\|^2\right]\right)^{1/2} \leq K\Delta^{(q+2)/2}
  \end{equation*}
  holds for the modelling error.
  
  \subsubsection*{Step 2: Itô isometry}
  Inserting the integral equation of \eqref{Omega SDE} and its truncated version, we get
  \begin{align*}
      &\left(\mathbb{E}\left[\|\Omega_{\Delta}-\Omega_{\Delta}^q\|^2\right]\right)^{1/2} \\
      &= \Bigg(\mathbb{E}\bigg[\Big\Vert\int_0^{\Delta}\sum_{k=q+1}^{\infty}\frac{B_k}{k!}\ad_{-\Omega_s}^k(K_s-\frac{1}{2}V_s^2)ds + \int_0^{\Delta}\sum_{k=q+1}^{\infty} \frac{B_k}{k!}\ad_{-\Omega_s}^k(V_s)dW_s\Big\Vert^2\bigg]\Bigg)^{1/2} \\
      &\leq \Bigg(\mathbb{E}\bigg[\Big\Vert\int_0^{\Delta}\sum_{k=q+1}^{\infty}\frac{B_k}{k!}\ad_{-\Omega_s}^k(K_s-\frac{1}{2}V_s^2)ds\Big\Vert^2\bigg]\Bigg)^{1/2} + \Bigg(\mathbb{E}\bigg[\Big\Vert\int_0^{\Delta}\sum_{k=q+1}^{\infty} \frac{B_k}{k!}\ad_{-\Omega_s}^k(V_s)dW_s\Big\Vert^2\bigg]\Bigg)^{1/2} \\
      &\leq \Bigg(\int_0^{\Delta}\mathbb{E}\bigg[\Big\Vert\sum_{k=q+1}^{\infty}\frac{B_k}{k!}\ad_{-\Omega_s}^k(K_s-\frac{1}{2}V_s^2)\Big\Vert^2\bigg]ds\Bigg)^{1/2} + \Bigg(\int_0^{\Delta}\mathbb{E}\bigg[\Big\Vert\sum_{k=q+1}^{\infty} \frac{B_k}{k!}\ad_{-\Omega_s}^k(V_s)\Big\Vert^2\bigg]ds\Bigg)^{1/2} \\
      &\leq \Bigg(\int_0^{\Delta}\mathbb{E}\bigg[\Big(\sum_{k=q+1}^{\infty}\frac{|B_k|}{k!}\big\|\ad_{-\Omega_s}^k(K_s-\frac{1}{2}V_s^2)\big\|\Big)^2\bigg]ds\Bigg)^{1/2} + \Bigg(\int_0^{\Delta}\mathbb{E}\bigg[\Big(\sum_{k=q+1}^{\infty} \frac{|B_k|}{k!}\big\|\ad_{-\Omega_s}^k(V_s)\big\|\Big)^2\bigg]ds\Bigg)^{1/2},
  \end{align*}
  where we also used the Minkowski inequality, the Itô isometry and the properties of a matrix norm. Now, the summands in the last line differ only in the input matrix of the adjoint operator.
  
  \subsubsection*{Step 3: Adjoint operator}
  We estimate the Frobenius norm of the adjoint operator of $V_s$ for a fixed $s\in [0,\Delta]$ and keep in mind that analogous estimates hold for the adjoint operator of $K_s-\frac{1}{2}V_s^2$. Since the Frobenius norm is submultiplicative, we have
  \begin{equation*}
      \|\ad_{-\Omega_s}(V_s)\| = \|[-\Omega_s,V_s]\| \leq \|\Omega_sV_s\|+\|V_s\Omega_s\|\leq 2\|\Omega_s\|\|V_s\|.
  \end{equation*}
  As a direct consequence, it holds
  \begin{equation*}
      \|\ad_{-\Omega_s}^k(V_s)\| \leq 2^k\|\Omega_s\|^k\|V_s\|,
  \end{equation*}
  which can also be shown via induction. Inserting this result in the expected value considered in the last line of the previous step, we get
  \begin{align*}
      \mathbb{E}\bigg[\Big(\sum_{k=q+1}^{\infty}\frac{|B_k|}{k!}\big\|\ad_{-\Omega_s}^k(V_s)\big\|\Big)^2\bigg] 
      &\leq \mathbb{E}\bigg[\Big(\sum_{k=q+1}^{\infty}\frac{|B_k|}{k!}2^k\|\Omega_s\|^k\|V_s\|\Big)^2\bigg]\\
      &= \|V_s\|^2\,\mathbb{E}\bigg[\Big(\sum_{k=q+1}^{\infty}\frac{|B_k|}{k!}2^k\|\Omega_s\|^k\Big)^2\bigg].
  \end{align*}
  
  \subsubsection*{Step 4: Estimate for the remainder}
  It is known that the Bernoulli numbers are implicitly defined by $\sum_{k=0}^{\infty}(B_k/k!)x^k=x/(e^x-1)$. Inserting the absolute values of the Bernoulli numbers instead, it holds
  \begin{equation*}
      \sum_{k=0}^{\infty}\frac{|B_k|}{k!}x^k = \frac{x}{2}\left(1+\cot\left(\frac{x}{2}\right)\right)+2.
  \end{equation*}
  Let $f:I\to \mathbb{R}$, $x \mapsto \frac{x}{2}\left(1+\cot\left(\frac{x}{2}\right)\right)+2$ with $I=\{x \in \mathbb{R}: \frac{x}{2\pi}\not\in\mathbb{Z}\}$. Applying Taylor's theorem to the function $f$ at the point 0 reads
  \begin{equation*}
      f(x)= \sum_{k=0}^{q}\frac{f^{(k)}(0)}{k!}x^k + R_q(x), \quad R_q(x)=\frac{f^{(q+1)}(\xi)}{(q+1)!}x^{q+1},
  \end{equation*}
  where we consider the Lagrange form of the remainder for some real number $\xi$ between 0 and $x$.
  
  Setting $x=2\|\Omega_s\|$ and recalling that the expression \eqref{eq: dexpinv} only converges for $\|\Omega\|<\pi$, we now consider $\left.f\right|_{\Tilde{I}}:\Tilde{I}\to \mathbb{R}$, $x \mapsto \frac{x}{2}\left(1+\cot\left(\frac{x}{2}\right)\right)+2$ with $\Tilde{I}=\{x \in \mathbb{R}: |x|<\pi\}$. The restriction of $f$ to $\Tilde{I}$ is bounded, in particular there exists an upper bound $M_q$ such that $|\left.f\right|_{\Tilde{I}}^{(q+1)}(\xi)|\le M_q$ for all $\xi$ between 0 and $x$.
  Moreover, the following estimate for the remainder holds 
  \begin{equation*}
      |R_q(x)|=\left|\frac{\left.f\right|_{\Tilde{I}}^{(q+1)}(\xi)}{(q+1)!}x^{q+1}\right|\leq \frac{M_q}{(q+1)!}|x|^{q+1} \leq \frac{M_q}{(q+1)!}(2\|\Omega_s\|)^{q+1}.
  \end{equation*}
  
  Using this estimate in the expected value of the last line of the previous step results in
  \begin{align*}
      \mathbb{E}\bigg[\Big(\sum_{k=q+1}^{\infty}\frac{|B_k|}{k!}(2\|\Omega_s\|)^k\Big)^2\bigg]&\leq \mathbb{E}\bigg[\Big(\frac{M_q}{(q+1)!}(2\|\Omega_s\|)^{q+1}\Big)^2\bigg]\\
      &= \left(\frac{2^{q+1}M_q}{(q+1)!}\right)^2 \mathbb{E}\left[\|\Omega_s\|^{2q+2}\right].
  \end{align*}
  
  \subsubsection*{Step 5: Itô-Taylor expansion}
  The goal of this step is to find an estimate for $\mathbb{E}\left[\|\Omega_s\|^{2q+2}\right]$. For this purpose, we examine the following Itô-Taylor expansion
  \begin{equation*}
      \Omega_{\Delta}=\Gamma(\Omega_0)\int_0^\Delta dW_s + R_{\Delta} = V_0W_{\Delta}+R_{\Delta},\quad \mathbb{E}\left[\|R_{\Delta}\|^2\right] \leq C_1\Delta^2,
  \end{equation*}
  where $C_1$ is a finite constant, for details see \cite[Proposition~5.9.1]{KlPl92}.
  Hence, the Frobenius norm of $\Omega_s$ can be estimated by
  \begin{equation*}
      \|\Omega_s\|=\|V_0W_s+R_s\|\leq \|V_0\||W_s|+\|R_s\|.
  \end{equation*}
  This result allows us to use the formula for the moments of the Wiener and the estimate for the remainder of the Itô-Taylor expansion,
  \begin{align*}
      \mathbb{E}\left[\|\Omega_s\|^{2q+2}\right] 
      &\leq 2^{2q+1}\left(\|V_0\|^{2q+2}\mathbb{E}\left[W_s^{2(q+1)}\right]+\mathbb{E}\left[\|R_s\|^{2(q+1)}\right]\right)\\
      &\leq 2^{2q+1}\left(\|V_0\|^{2q+2}\frac{(2(q+1))!}{2^{q+1}(q+1)!}s^{q+1}+C_1s^{2(q+1)}\right).
  \end{align*}
  
  \subsubsection*{Step 6: Overall estimate}
  Gathering the results of the previous steps and inserting a Taylor expansion for $V_s$ where $C_2<\infty$ reads
  \begin{align*}
      &\mathbb{E}\bigg[\Big(\sum_{k=q+1}^{\infty} \frac{|B_k|}{k!}\big\|\ad_{-\Omega_s}^k(V_s)\big\|\Big)^2\bigg]\\
      &\leq \|V_s\|^2\left(\frac{2^{q+1}M_q}{(q+1)!}\right)^2 2^{2q+1}\left(\|V_0\|^{2q+2}\frac{(2(q+1))!}{2^{q+1}(q+1)!}s^{q+1}+C_1s^{2(q+1)}\right)\\
      &\leq \left(\|V_0\|+ C_2s\right)^2\left(\frac{2^{q+1}M_q}{(q+1)!}\right)^2 2^{2q+1}\left(\|V_0\|^{2q+2}\frac{(2(q+1))!}{2^{q+1}(q+1)!}s^{q+1}+C_1s^{2(q+1)}\right)\\
      &=\mathcal{O}(s^{q+1}).
  \end{align*}
  Thus, it holds
  \begin{equation*}
      \Bigg(\int_0^{\Delta}\mathbb{E}\bigg[\Big(\sum_{k=q+1}^{\infty} \frac{|B_k|}{k!}\big\|\ad_{-\Omega_s}^k(V_s)\big\|\Big)^2\bigg]ds\Bigg)^{1/2} =\mathcal{O}(\Delta^{(q+2)/2}).
  \end{equation*}
  Analogously, one can show that
  \begin{equation*}
      \Bigg(\int_0^{\Delta}\mathbb{E}\bigg[\Big(\sum_{k=q+1}^{\infty} \frac{|B_k|}{k!}\big\|\ad_{-\Omega_s}^k(K_s-\frac{1}{2}V_s^2)\big\|\Big)^2\bigg]ds\Bigg)^{1/2} =\mathcal{O}(\Delta^{(q+2)/2}),
  \end{equation*}
  which concludes the proof.\qed
\end{proof}

Note that due to the definition of the Cayley map as a finite product of matrices no such theorem is needed if $\psi=\cay$ is chosen as the local parametrization in Algorithm~\ref{algo}.

We further point out that Theorem~\ref{theorem: index} is in accordance with our results of Section~\ref{subsec: order1}, where the geometric Euler-Maruyama scheme~\eqref{eq: gEM exp} can be interpreted as a stochastic RKMK method with $\gamma=1$ and $q=0$.

\section{Numerical examples}\label{sec4}
In the following we provide numerical examples which illustrate the effectiveness of the proposed geometric methods, 
firstly, by simulating the strong convergence order of the proposed schemes and secondly, by showing the Lie group structure preservation of our methods.

For checking the convergence order, we set $G=$ SO(3) and $\mathfrak{g}=\mathfrak{so}(3)$. 
In order to ensure the conditions of Theorem~\ref{theorem: so(n)} we have used the set up of matrices $K_t$ and $V_t$ proposed by Muniz et al.~\cite{Muniz20}. 
Specifically, we chose the time-dependent functions 
\begin{equation*}
    f_1(t) = \cos(t), \quad
    f_2(t) = \sin(t), \quad
    f_3(t) = 1 + t + t^2 + t^3,
\end{equation*}
to compute a skew-symmetric matrix $V_t$ as a linear combination,
\begin{equation*}
    V_t = f_1(t) G_1 + f_2(t)G_2 + f_3(t)G_3,
\end{equation*}
where $G_i$, $i=1,2,3$ are the following generators of the Lie algebra $\mathfrak{so}(3)$,
\begin{equation*}
    G_1 = 
    \begin{pmatrix}
    0 & -1 & 0 \\ 1 & 0 & 0 \\ 0 & 0 & 0
    \end{pmatrix}, \quad 
    G_2 =
    \begin{pmatrix}
    0 & 0 & -1 \\ 0 & 0 & 0 \\ 1 & 0 & 0
    \end{pmatrix}, \quad
    G_3 =
    \begin{pmatrix}
    0 & 0 & 0 \\ 0 & 0 & -1 \\ 0 & 1 & 0
    \end{pmatrix}.
\end{equation*}
Note that the functions $f_i$, $i=1,2,3$ can be chosen arbitrarily.
We then set the matrix $K_t$ as the lower triangular matrix of $V_t^2$ where the diagonal entries of $K_t$ are 0.5 times the diagonal entries of $V_t^2$ such that $K_t + K_t^{\top} = V_t^2$.

We simulated $M=1000$ different paths of two independent realizations of a standard normally distributed random variable, $U_1, U_2 \sim\mathcal{N}(0,1)$.
Then, the random variables used in the numerical method step in Algorithm~\ref{algo} were simulated as $\widehat{\Delta W} = U_1\sqrt{\Delta}$ and $\widehat{\Delta Z} = \frac{1}{2}\Delta(\widehat{\Delta W} + U_2\sqrt{\frac{\Delta}{3}})$. 
The absolute error as defined in~\eqref{eq: order} was estimated by using the Frobenius norm at $t_j=T$, i.e.\ by 
\begin{equation*}
    \frac{1}{M}\sum_{i=1}^{M} \Bigl(\bigl\|Q_{T,i}^{\text{ref}}-Q_{T,i}^{\Delta}\bigr\|_F\Bigr)
\end{equation*}
where the approximations $Q_{T,i}^{\Delta}$ were obtained by using Algorithm~\ref{algo} with step sizes $\Delta = 2^{-14}, 2^{-13}, 2^{-12}, 2^{-11}, 2^{-10}, 2^{-9}$ and for the reference solution $Q_{T,i}^{\text{ref}}$ we used the same method with $\psi=\cay$ and step size $\Delta=2^{-16}$, respectively. 

\begin{figure}
\hspace{-0.8cm}
    \includegraphics[scale=0.58]{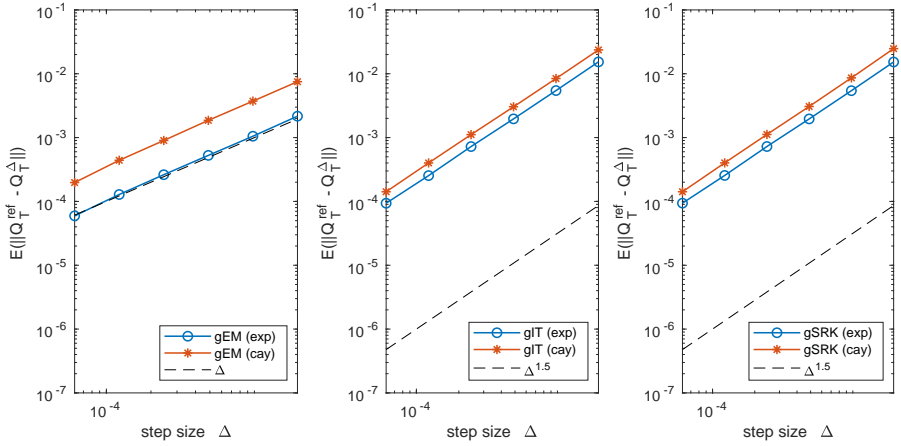}
    \caption{Simulation of the strong convergence order for $M=1000$ paths. 
    Left: Geometric Euler-Maruyama scheme (gEM). Center: Geometric version of the Itô-Taylor scheme of strong order 1.5 (gIT). 
    Right: Geometric version of Rößler's stochastic Runge-Kutta scheme of strong order 1.5 (gSRK).}
    \label{plot: conv}
\end{figure}

A log-log plot of the estimation of the absolute error against the step sizes can be viewed in Figure~\ref{plot: conv}. It indicates the strong order of convergence claimed in the sections above for the geometric Euler-Maruyama scheme~\eqref{eq: gEM}, the geometric version of the Itô-Taylor scheme~\eqref{eq: ItôTaylor1.5} and the geometric stochastic Runge-Kutta scheme~\eqref{eq: Roessler1.5}.

Examples from financial mathematics and multibody system dynamics verify that the  structure-preserving methods derived above can be applied in practice. 

In the first example we apply our methods of strong order $\gamma=1.5$ to an SDE on SO(2) in the context of stochastic correlation modelling.  
The second example shows how our methods can be used in the modeling of rigid bodies, e.g.\ satellites. 
Although, we have restricted our research for this paper to considering only linear SDEs on Lie groups, the second example shows that our methods can also be applied to nonlinear SDEs on e.g.\ SO(3).

\subsection{A stochastic correlation model}
Let us assume that a risk manager retrieves from the middle office's reporting system an initial value of the correlation between two assets and a density function of the considered correlation. Moreover, we assume that the risk manager was given the task to generate correlation matrices that not only approximate the given density function but also respect the stochastic behaviour of correlations.

This problem can be solved by the stochastic correlation model presented in \cite{Muniz20}. The main ideas of the approach are outlined in the following.

\begin{figure}
    \centering
   \includegraphics[scale=0.45]{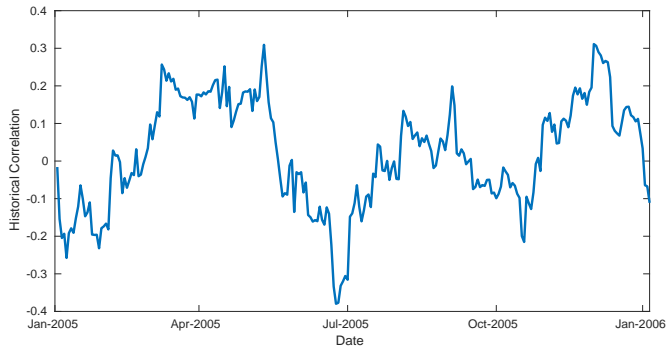}
   \caption{The 30-day historical correlations between S\&P 500 and Euro/US-Dollar exchange rate, source of data: 
   \texttt{www.yahoo.com}.}
   \label{fig:histcorr}
\end{figure}

For this example we consider historical prices of the S\&P 500 index and the Euro/US-Dollar exchange rate and compute moving correlations with a window size of 30 days to obtain correlations from January 03, 2005 to January 06, 2006 (see Figure~\ref{fig:histcorr}).

\begin{figure}
\centering
   \includegraphics[scale=0.5]{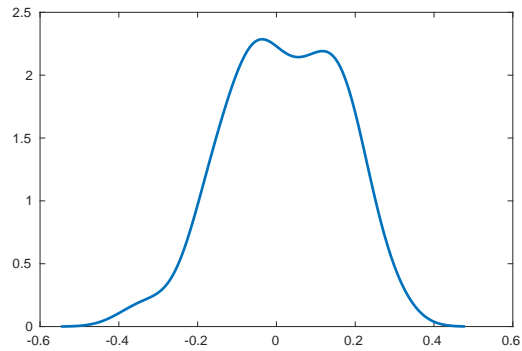}
   \caption{Empirical density function of the historical correlation between S\&P 500 and Euro/US-Dollar exchange rate, computed with the MATLAB function \texttt{ksdensity}.}
   \label{fig:histdensity}
\end{figure}

The corresponding initial correlation matrix calculated from this data and imputed to the risk manager is
\begin{equation*}
	R_0^{\rm hist}=
	\begin{pmatrix}
	1 & -0.0159 \\
	-0.0159 & 1 
	\end{pmatrix}.
\end{equation*}
Furthermore, we estimate a density function from the historical data using kernel smoothing functions, which is also plotted in Figure~\ref{fig:histdensity}. 
For more details on the density estimation see \cite{BoAz97}. 

As a first step, we focus on covariance matrices $P_t$, $t\ge 0$. 
The authors of \cite{Teng19} utilised the principal axis theorem and defined the \textit{covariance flow}
\begin{equation} \label{eq: covflow}
    P_t = Q_t^{\top}P_0Q_t, \quad t\ge 0,
\end{equation}
where $P_0$ is the initial covariance matrix computed based on $R_0^{\rm hist}$ and $Q_t$ is an orthogonal matrix which without loss of generality can be assumed to have determinant +1, i.e.\ $Q_t \in \mathrm{SO}(2)$. 
Following the approach in \cite{Muniz20} the matrix $Q_t$ is now assumed to be driven by the SDE~\eqref{Q SDE} which can be solved by using Algorithm~\ref{algo}. 
With the resulting matrices approximations of $P_t$ can be computed with~\eqref{eq: covflow}, which can then be transformed to corresponding correlation matrices
\begin{equation*}
    R_t=\Sigma_t^{-1} P_t \Sigma_t^{-1}
\end{equation*}
with $\Sigma_t=\bigl(\diag(P_t)\bigr)^{\frac{1}{2}}$.

At last, a density function is estimated from this \textit{correlation flow} and the free parameters involved are calibrated such that the density function matches the density function from the historical data, see \cite{Muniz20} for details.

\begin{figure}
\centering
   \includegraphics[scale=0.5]{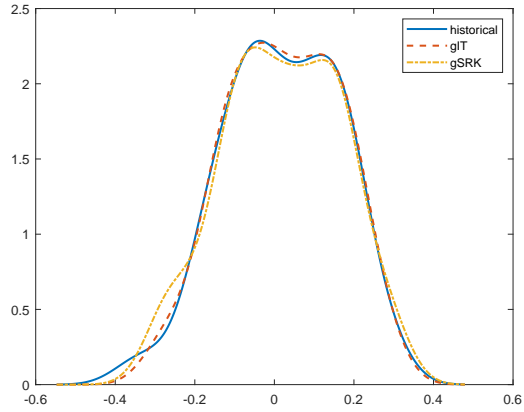}
   \caption{Empirical density function of the historical correlation and the correlation flow between S\&P 500 and Euro/US-Dollar exchange rate.}
   \label{fig:flowdensity}
\end{figure}

We executed this procedure using the geometric Itô-Taylor scheme~\eqref{eq: ItôTaylor1.5} with $\psi=\cay$ (gIT) and the geometric R\"oßler scheme~\eqref{eq: Roessler1.5} with $\psi=\exp$ and truncation index $q=1$ (gSRK) in the second step of Algorithm~\ref{algo}, respectively. 
The results are plotted in Figure~\ref{fig:flowdensity}, which shows that both density functions approximate the density function of the historical data quite well.

\subsection{The stochastic rigid body problem}
Consider a free rigid body, whose centre of mass is at the origin. Let the vector $y=(y_1,y_2,y_3)^{\top}$ represent the angular momentum in the body frame and $I_1$, $I_2$ and $I_3$ be the principal moments of inertia \cite{Mars99}. 
Then the motion of this free rigid body is described by the Euler equations
\begin{equation*}
    \Dot{y} = V(y)y, \quad V(y) = 
    \begin{pmatrix}
    0 & y_3/I_3 & -y_2/I_2 \\
    -y_3/I_3 & 0 & y_1/I_1 \\
    y_2/I_2 & -y_1/I_1 & 0
    \end{pmatrix}.
\end{equation*}
We suppose that the rigid body is perturbed by a Wiener process $W_t$ and compute a matrix $K(y)$ such that the dynamics are kept on the manifold, i.e.\ we compute $K(y)$ from the condition $K(y) + K^{\top}(y) = V^2(y)$. 
Consequently, we regard the Itô SDE
\begin{equation}\label{eq: rb}
    dy = K(y)y\, dt + V(y)y\, dW_t,
\end{equation}
where the solution evolves on the unit sphere if the initial value $y_0$ satisfies $|y_0|=1$.
Note that stochastic versions of the rigid body problem have already been considered in \cite{MaWi08} and \cite{Wang20} but as Stratonovich SDEs.

Since the solution of~\eqref{eq: rb} can also be written as $y=Qy_0$ where $Q\in \mathrm{SO}(3)$, we focus on the nonlinear matrix SDE
\begin{equation}\label{eq: rb so3}
    dQ = K(Q)Q\, dt + V(Q)Q\, dW_t, \quad Q_0 = I_{3\times 3}.
\end{equation}
The coefficients of the corresponding SDE in the Lie algebra~\eqref{Omega SDE} read
\begin{equation}\label{eq: coeff}
    A(\Omega)=d\psi^{-1}_{\Omega}\Bigl(K\bigl(\psi(\Omega)Q_0\bigr) -\frac{1}{2}V^2\bigl(\psi(\Omega)Q_0\bigr)\Bigr), \quad
    \Gamma(\Omega)= d\psi^{-1}_{\Omega}\Bigl(V\bigl(\psi(\Omega)Q_0\bigr)\Bigr).
\end{equation}
Now, SDE~\eqref{eq: rb} can be solved by applying Algorithm~\ref{algo} to the SDE~\eqref{eq: rb so3}. 
Note that we deal with right multiplication of the solution $Q$ on the right hand side of~\eqref{eq: rb so3} instead of left multiplication as in~\eqref{Q SDE}. 
As a consequence, the sign of the index of the operator $d\psi^{-1}$ is changed in \eqref{eq: coeff} and the solution of the Projection step in Algorithm~\ref{algo} should be $Q_{j+1}=\psi(\Omega_{j+1})Q_j$.
We refer to~\cite{Munt98} for more details on this matter.

\begin{figure}
\centering
    \includegraphics[scale=0.7]{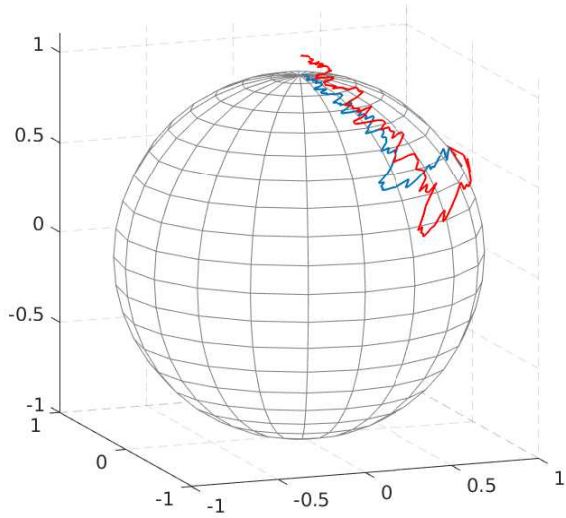}
    \caption{Sample paths of the geometric Euler-Maruyama (blue) and the traditional Euler-Maruyama scheme (red) applied to \eqref{eq: rb}.}
    \label{plot: rigid body}
\end{figure}

\begin{figure}
\centering
    \includegraphics[scale=0.6]{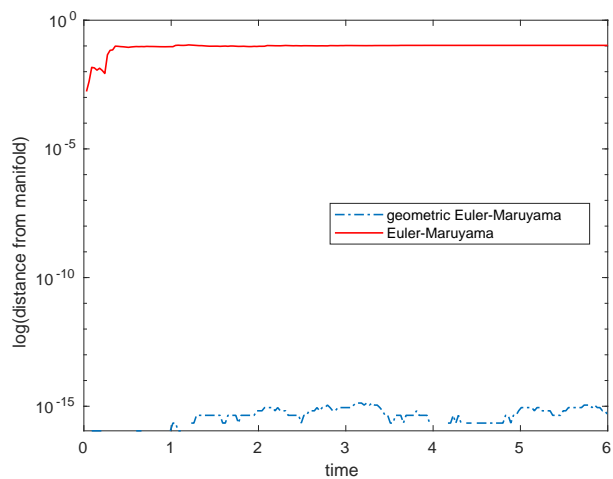}
    \caption{Log-distance of the numerical solutions to the unit sphere.}
    \label{plot: distance}
\end{figure}

In Figure~\ref{plot: rigid body} we simulated 200 steps of the trajectory of \eqref{eq: rb} with a step size of $\Delta=0.03$ by using Algorithm~\ref{algo} with the initial values $y_0=(\sin(1.1),0,\cos(1.1))^{\top}$ and the moments of inertia $I_1=2$, $I_2=1$ and $I_3=2/3$. 
For the numerical method step of Algorithm~\ref{algo} we used the Euler-Maruyama scheme with $\psi=\cay$. 
Emphasizing the structure-preserving character of Algorithm~\ref{algo} we also plotted a sample path of the traditional Euler-Maruyama scheme applied directly to \eqref{eq: rb}, whose trajectory clearly fails to stay on the manifold. 
This phenomenon can also be viewed in Figure~\ref{plot: distance} where we visualize the distance of the approximate solutions from the manifold.

\section{Conclusion}\label{sec5}

We have presented stochastic Lie group methods for linear Itô SDEs on matrix Lie groups that have a higher strong convergence order than the known \textit{geometric Euler-Maruyama} scheme. 
Based on RKMK methods for ODEs on Lie groups, we have proven a condition on the truncation index of the inverse of $d\exp(H)$ such that the stochastic RKMK method inherits the convergence order of the underlying SRK. 
Additionally, we have shown examples for the application of our methods in mechanical engineering and in financial mathematics.

Our methods require further investigations for the application to nonlinear Itô SDEs on matrix Lie groups, 
which we consider as future work. Moreover, we have restricted our research for this paper to the strong convergence order.
In future research, an investigation on the weak convergence order of stochastic Lie group methods will also be conducted.

\begin{acknowledgements}
The authors would like to thank Martin Friesen (Dublin City University) for the in-depth discussions that improved the content of this paper.

The work of the authors was partially supported by the bilateral German-Slovakian Project 
\textit{MATTHIAS -- Modelling and Approximation Tools and Techniques for Hamilton-Jacobi-Bellman equations in finance and Innovative Approach to their Solution}, financed by DAAD and the Slovakian
Ministry of Education. Further the authors acknowledge partial support from the bilateral German-Portuguese Project
\textit{FRACTAL -- FRActional models and CompuTationAL Finance} financed by DAAD and the CRUP - Conselho de Reitores
das Universidades Portuguesas.
\end{acknowledgements}



\appendix
\section{Proofs}

\begin{theorem}
The solution of~\eqref{Q SDE} can be written as $Q_t=Q_0 \,\psi(\Omega_t)$, where $\Omega_t$ obeys the SDE~\eqref{Omega SDE} with coefficients given by
\begin{equation*}
   A(\Omega_t) = d\psi_{-\Omega_t}^{-1}\Bigl(K_t - \frac{1}{2}V_t^2 - \frac{1}{2}C(\Omega_t)\Bigr),\quad \Gamma(\Omega_t) = d\psi_{-\Omega_t}^{-1}(V_t)
\end{equation*}
with 
\begin{equation}\label{eq:coeff}
   C(\Omega_t) = \Bigl(\frac{d}{d\Omega_t}d\psi_{-\Omega_t}\bigl(\Gamma(\Omega_t)\bigr)\Bigr)\Gamma(\Omega_t).
\end{equation}
\end{theorem}

\begin{proof}
Let $f\colon\mathfrak{g}\to G$ be the chosen local parametrization 
such that $Q_t=f(\Omega_t)$ with $f(\Omega_t)=Q_0\psi(\Omega_t)$. 
Due to the Itô rules $dQ_t=d(f(\Omega_t))$ is given exactly by the first two terms of the Taylor expansion
\begin{equation*}
\begin{split}
    dQ_t &= \frac{d}{d\epsilon}f(\Omega_t + \epsilon d\Omega_t)\Bigr|_{\epsilon=0} + \frac{1}{2} \frac{d^2}{d\epsilon^2}f(\Omega_t + \epsilon d\Omega_t)\Bigr|_{\epsilon=0} \\
    &= \Bigl(\frac{d}{d\Omega_t}f(\Omega_t)\Bigr)\,d\Omega_t + \frac{1}{2}\Bigl(\frac{d^2}{d\Omega_t^2}f(\Omega_t)\Bigr)\,(d\Omega_t)^2 \\
    &= \left(\Bigl(\frac{d}{d\Omega_t}f(\Omega_t)\Bigr)A(\Omega_t) + \frac{1}{2}\Bigl(\frac{d^2}{d\Omega_t^2}f(\Omega_t)\Bigr)\Gamma^2(\Omega_t)\right) dt + \Bigl(\frac{d}{d\Omega_t}f(\Omega_t)\Bigr)
    \Gamma(\Omega_t)\, dW_t
\end{split}
\end{equation*}
where we have used the fact that 
\begin{equation*}
    (d\Omega_t)^2=\bigl(A(\Omega_t)\,dt+\Gamma(\Omega_t)\,dW_t\bigr)^2= \Gamma^2(\Omega_t)\,dt.
\end{equation*}
For both $\psi=\exp$ and $\psi=\cay$ it holds that
\begin{equation*}
\Bigl(\frac{d}{d\Omega} \psi(\Omega)\Bigr)H = \bigl(d\psi_{\Omega}(H)\bigr)\psi(\Omega)=\psi(\Omega)\left(d\psi_{-\Omega}(H)\right),
\end{equation*}
which we use to specify the first part of the drift coefficient
\begin{equation*}
 \begin{split}
    \Bigl(\frac{d}{d\Omega_t}f(\Omega_t)\Bigr)A(\Omega_t) 
    &= Q_0 \,\Bigl(\frac{d}{d\Omega_t}\psi(\Omega_t)\Bigr)A(\Omega_t) 
    = Q_0 \,\psi(\Omega_t)d\psi_{-\Omega_t}\bigl(A(\Omega_t)\bigr) \\
    &= Q_t\,d\psi_{-\Omega_t}\bigl(A(\Omega_t)\bigr).
\end{split}
\end{equation*}
Analogously, we have $\bigl(\frac{d}{d\Omega_t}f(\Omega_t)\bigr)\,\Gamma(\Omega_t)= Q_t d\psi_{-\Omega_t}\bigl(\Gamma(\Omega_t)\bigr)$. 
For the second derivative we obtain
\begin{equation*}
\begin{split}
    \Bigl(\frac{d^2}{d\Omega_t^2}f(\Omega_t)\Bigr)\Gamma^2(\Omega_t) 
    &= \Bigl( \frac{d}{d\Omega}Q_t d\psi_{-\Omega_t}\bigl(\Gamma(\Omega_t)\bigr)\Bigr)\Gamma(\Omega_t) \\
    &= \left(\Bigl(\frac{d}{d\Omega_t}Q_t\Bigr)\Gamma(\Omega_t)\right)d\psi_{-\Omega_t}\bigl(\Gamma(\Omega_t)\bigr) \\
    &\qquad+ Q_t\Bigl(\frac{d}{d\Omega_t}d\psi_{-\Omega_t}\bigl(\Gamma(\Omega_t)\bigr)\Bigr)\Gamma(\Omega_t) \\
    &= Q_t\,\Bigl(d\psi_{-\Omega_t}\bigl(\Gamma(\Omega_t)\bigr)\Bigr)^2 + Q_t C(\Omega_t),
\end{split}
\end{equation*}
where $C(\Omega_t)=\Bigl(\frac{d}{d\Omega_t}d\psi_{-\Omega_t}\bigl(\Gamma(\Omega_t)\bigr)\Bigr)\Gamma(\Omega_t)$. 
Comparing these results with SDE~\eqref{Q SDE} we get
\begin{equation*}
    V_t = d\psi_{-\Omega_t}\bigl(\Gamma(\Omega_t)\bigr) \quad \text{and} \quad K_t = d\psi_{-\Omega_t}\bigl(A(\Omega_t)\bigr) + \frac{1}{2}\left(d\psi_{-\Omega_t}\bigl(\Gamma(\Omega_t)\bigr)\right)^2 + \frac{1}{2}C(\Omega_t)
\end{equation*}
and thus
\begin{equation*}
    \Gamma(\Omega_t) = d\psi_{-\Omega_t}^{-1}(V_t) \quad \text{and} \quad A(\Omega_t) = d\psi_{-\Omega_t}^{-1}\Bigl(K_t - \frac{1}{2}V_t^2 - \frac{1}{2}C(\Omega_t)\Bigr).
\end{equation*}\qed
\end{proof}

\begin{lemma}
For $\psi=\cay$ the coefficient~\eqref{eq:coeff} is given by
\begin{equation}
    C(\Omega_t) = V_t\Omega_t V_t. 
\end{equation}
\end{lemma}

\begin{proof}
\begin{equation*}
\begin{split}
&\left(\frac{d}{d\Omega_t}d\cay_{-\Omega_t}(H)\right)\Tilde{H}\\
   &= \left(2 \frac{d}{d\Omega_t}\Bigl((I+\Omega_t)^{-1}H(I-\Omega_t)^{-1}\Bigr)\right)\Tilde{H}\\
   &= 2\left(\Bigl(\frac{d}{d\Omega_t}(I+\Omega_t)^{-1}\Bigr)\Tilde{H}\right)H(I-\Omega_t)^{-1} + 2(I+\Omega_t)^{-1}H\left(\frac{d}{d\Omega_t}(I-\Omega_t)^{-1}\right)\Tilde{H}\\
   &= -2(I+\Omega_t)^{-1}\Tilde{H}(I+\Omega_t)^{-1}H(I-\Omega_t)^{-1} + 2(I+\Omega_t)^{-1}H(I-\Omega_t)^{-1}\Tilde{H}(I-\Omega_t)^{-1}
\end{split}
\end{equation*}
Inserting for both $H$ and $\Tilde{H}$ the diffusion coefficient 
$\Gamma(\Omega_t)=d\cay_{-\Omega_t}^{-1}(V_t) = \frac{1}{2}(I+\Omega_t)V_t(I-\Omega_t)$ we get
\begin{equation*}
\begin{split}
\biggl(&\frac{d}{d\Omega_t}d\cay_{-\Omega_t}\bigl(\Gamma(\Omega_t)\bigr)\biggr)\Gamma(\Omega_t)\\
   &= -2(I+\Omega_t)^{-1}\Gamma(\Omega_t)(I+\Omega_t)^{-1}\Gamma(\Omega_t)(I-\Omega_t)^{-1}\\
   &\qquad + 2(I+\Omega_t)^{-1}\Gamma(\Omega_t)(I-\Omega_t)^{-1}\Gamma(\Omega_t)(I-\Omega_t)^{-1} \\
   &= -\frac{1}{2}(I+\Omega_t)^{-1}(I+\Omega_t)V_t(I-\Omega_t)(I+\Omega_t)^{-1}(I+\Omega_t)V_t(I-\Omega_t)(I-\Omega_t)^{-1} \\
   & \qquad+\frac{1}{2}(I+\Omega_t)^{-1}(I+\Omega_t)V_t(I-\Omega_t)(I-\Omega_t)^{-1}(I+\Omega_t)V_t(I-\Omega_t)(I-\Omega_t)^{-1} \\
   &= -\frac{1}{2}V_t(I-\Omega_t)V_t +\frac{1}{2}V_t(I+\Omega_t)V_t \\
   &= V_t\Omega_t V_t.
\end{split}
\end{equation*}\qed
\end{proof}

\section{Matrix derivatives}
In this section we provide the matrix derivatives that we used in the geometric version of the Itô-Taylor scheme of strong order $\gamma=1.5$ (see \eqref{eq: ItôTaylor1.5}).

\subsection{Derivatives for $\psi=\cay$}\label{subsec: cay derivatives}
Computing the derivative of \eqref{eq: dcayinv} in the direction of an arbitrary matrix $\Tilde{H}$, we get
\begin{equation}\label{eq: ddcay}
\begin{split}
    \left(\frac{d}{d\Omega}d\cay_{-\Omega}^{-1}(H)\right)\Tilde{H} 
    &= \left(\frac{1}{2}\frac{d}{d\Omega}(H-H\Omega+\Omega H - \Omega H\Omega)\right)\Tilde{H} \\
    &= \frac{1}{2} \frac{d}{dt}\bigl(H - H(\Omega+t\Tilde{H}) + (\Omega+t\Tilde{H})H \\
    &\qquad\left. - (\Omega+t\Tilde{H})H(\Omega+t\Tilde{H})\bigr)\right|_{t=0} \\
    &= \frac{1}{2}(-H\Tilde{H} + \Tilde{H}H - \Omega H \Tilde{H} - \Tilde{H}H\Omega).
    \end{split}
\end{equation}

Subsequently, the second directional derivative reads
\begin{equation}\label{eq: dddcay}
    \left(\frac{d^2}{d\Omega^2}d\cay_{-\Omega}^{-1}(H)\right)\Tilde{H}^2 
    = -\Tilde{H}H\Tilde{H}.
\end{equation}

Inserting $H=V$ and $\Tilde{H}=\Gamma(\Omega)=d\cay_{-\Omega}^{-1}(V)$, we obtain
\begin{align*}
    \left(\frac{d}{d\Omega}\Gamma(\Omega)\right)\Gamma(\Omega)&=\left(\frac{d}{d\Omega}d\cay_{-\Omega}^{-1}(V)\right)\Gamma(\Omega)\\
    &= \frac{1}{2}(-V\Gamma(\Omega) + \Gamma(\Omega)V - \Omega V\Gamma(\Omega) -\Gamma(\Omega)V\Omega)\\
    &= \frac{1}{2}(-V\Omega V + V \Omega V \Omega - \Omega V \Omega V + \Omega V \Omega V \Omega).
\end{align*}

Similar expressions are obtained for $\left(\frac{d}{d\Omega}A(\Omega)\right)\Gamma(\Omega)$, $\left(\frac{d}{d\Omega}A(\Omega)\right)A(\Omega)$ and $\left(\frac{d}{d\Omega}\Gamma(\Omega)\right)A(\Omega)$ by inserting $H=K-\frac{1}{2}V^2$ and $\Tilde{H}=\Gamma(\Omega)$, $H=K-\frac{1}{2}V^2$ and $\Tilde{H}=A(\Omega)$ and $H=V$ and $\Tilde{H}=A(\Omega)$ in \eqref{eq: ddcay}, respectively.

Proceed accordingly to compute the second derivatives $\left(\frac{d^2}{d\Omega^2}A(\Omega)\right)\Gamma^2(\Omega)$ and $\left(\frac{d^2}{d\Omega^2}\Gamma(\Omega)\right)\Gamma^2(\Omega)$.

\subsection{Derivatives for $\psi=\exp$}
In the following we present derivatives of \eqref{eq: dexp} up to $k=4$, i.e.\ of
\begin{align*}
    \sum_{k=0}^{4}\frac{B_k}{k!}\ad_{-\Omega}^k(H) 
	= \,&H - \frac{1}{2} \lbrack -\Omega, H \rbrack + \frac{1}{12} \bigl\lbrack - \Omega, \lbrack -\Omega,H \rbrack \bigr\rbrack - \frac{1}{720}\left[-\Omega,\Bigl\lbrack -\Omega \bigl\lbrack - \Omega, \lbrack -\Omega,H \rbrack \bigr\rbrack\Bigr\rbrack\right]\\
	= \,&H - \frac{1}{2}(H\Omega - \Omega H) + \frac{1}{12}(\Omega^2H+ H\Omega^2 - 2\Omega H\Omega) \\
	&- \frac{1}{720}(\Omega^4H - 4\Omega^3H\Omega + 6\Omega^2H\Omega^2 - 4\Omega H \Omega^3 + H\Omega^4).
\end{align*}

Computing the directional derivative we get
\begin{align*}
    \left(\frac{d}{d\Omega}\sum_{k=0}^{4}\frac{B_k}{k!}\ad_{-\Omega}^k(H)\right)\Tilde{H} =\,&-\frac{1}{2}(H\Tilde{H} - \Tilde{H}H) \\
    &+ \frac{1}{12}(\Omega\Tilde{H}H + \Tilde{H}\Omega H + H\Omega\Tilde{H} + H\Tilde{H}\Omega - 2\Tilde{H}H\Omega - 2\Omega H\Tilde{H}) \\
    &-\frac{1}{720}\bigl(\Tilde{H}\Omega^3H + \Omega\Tilde{H}\Omega^2H + \Omega^2\Tilde{H}\Omega H + \Omega^3\Tilde{H}H\\
    &-4(\Omega\Tilde{H}\Omega H \Omega + \Tilde{H}\Omega^2H\Omega + \Omega^2\Tilde{H}H\Omega + \Omega^3\Tilde{H}H)\\
    &+6(\Omega^2H\Omega\Tilde{H} + \Omega^2H\Tilde{H}\Omega + \Omega\Tilde{H}H\Omega^2 + \Tilde{H}\Omega H\Omega^2)\\
    &-4(\Omega H\Omega\Tilde{H}\Omega + \Omega H\Tilde{H}\Omega^2 + \Omega H\Omega^2\Tilde{H} + \Tilde{H}H\Omega^3)\\
    &+ H\Tilde{H}\Omega^3 + H\Omega\Tilde{H}\Omega^2 + H\Omega^2\Tilde{H}\Omega + H\Omega^3\Tilde{H}\bigr).
\end{align*}
Whereas the second directional derivative is given by 
\begin{align*}
    \biggl(\frac{d^2}{d\Omega^2}&\sum_{k=0}^{4}\frac{B_k}{k!}\ad_{-\Omega}^k(H)\biggr)\Tilde{H}^2\\ =\, &\frac{1}{6}(\Tilde{H}^2H + H\Tilde{H}^2 - 2\Tilde{H}H\Tilde{H}) \\
    &-\frac{1}{360}\bigl( (\Omega \Tilde{H} \Omega \Tilde{H} H + \Omega \Tilde{H}^2 \Omega H + \Tilde{H}^2 \Omega^2 H + \Omega^2 \Tilde{H}^2 H + \Tilde{H} \Omega \Tilde{H} \Omega H + \Tilde{H} \Omega^2\Tilde{H} H)\\
    &-4 (\Omega \Tilde{H} \Omega H \Tilde{H} + \Omega \Tilde{H}^2 H \Omega + \Tilde{H}^2 \Omega H \Omega + \Tilde{H} \Omega^2 H \Tilde{H} + \Tilde{H} \Omega H \Omega + \Omega^2 \Tilde{H} H \Tilde{H})\\
    &+6 (\Omega^2 H \Tilde{H}^2 + \Omega \Tilde{H} H \Omega \Tilde{H} + \Tilde{H} \Omega H \Omega \Tilde{H} + \Omega \Tilde{H} H \Tilde{H} \Omega + \Tilde{H} \Omega H \Tilde{H} \Omega + \Tilde{H}^2 H \Omega^2)\\
    &-4 (\Omega H \Omega \Tilde{H}^2 + \Omega H \Tilde{H}^2 \Omega + \Tilde{H} H \Omega \Tilde{H} \Omega + \Omega H \Tilde{H} \Omega \Tilde{H} + \Tilde{H} H \Tilde{H} \Omega^2 + \Tilde{H} H \Omega^2 \Tilde{H})\\
    &+ (H \Tilde{H} \Omega \Tilde{H} \Omega + H \Tilde{H}^2 \Omega^2 + H \Tilde{H}\Omega^2\Tilde{H} + H \Omega \Tilde{H} \Omega \Tilde{H} + H \Omega \Tilde{H}^2 \Omega + H \Omega^2 \Tilde{H}^2)\bigr).
\end{align*}

Note that evaluating the derivatives at $\Omega_0=0_{n\times n}$ causes many summands to become zero, which makes computing higher summands ($k>4$) unnecessary.

The needed derivatives for the geometric Itô-Taylor scheme~\eqref{eq: ItôTaylor1.5} can be computed from the formulas above by inserting correspondingly into $H$ and $\Tilde{H}$ (see \ref{subsec: cay derivatives} for instructions).

\end{document}